\documentclass[reqno]{amsart}
\usepackage{amsfonts,amsmath,amsthm,amssymb,mathrsfs}
\usepackage{color}
\usepackage{amsmath,amssymb,amsfonts,bm}
\usepackage{graphicx}
\usepackage{CJKutf8}
\usepackage{newunicodechar}
\usepackage{xcolor}
\numberwithin{equation}{section}
\usepackage[pdfstartview=FitH,colorlinks=true]{hyperref}
\hypersetup{urlcolor=red, linkcolor=blue,citecolor=red, CJKbookmarks=true}
\allowdisplaybreaks

\theoremstyle{plain}
\newtheorem{theorem}{Theorem}[section]
\newtheorem{remark}[theorem]{Remark}
\newtheorem{lemma}[theorem]{Lemma}
\newtheorem{proposition}[theorem]{Proposition}

\numberwithin{equation}{section}
\begin{document}

\title[Landau equation]
{Analytic smoothing effect of linear landau equation with soft potential}

\author[Hao-Guang Li   \& Chao-Jiang Xu]
{Hao-Guang Li and Chao-Jiang Xu}
\address{Hao-Guang Li,
\newline\indent
School of Mathematics and Statistics, South-Central University for Nationalities
\newline\indent
430074, Wuhan, China}
\email{lihaoguang@scuec.edu.cn}

\address{Chao-Jiang Xu,
\newline\indent
Department of Mathematics, Nanjing University of Aeronautics and Astronautics, Nanjing 211106, China
\newline\indent
Key Laboratory of Mathematical Modelling and High Performance Computing of Air Vehicles (NUAA), MIIT, Nanjing 210016, China
\newline\indent
Universit\'e de Rouen, CNRS UMR 6085, Laboratoire de Math\'ematiques
\newline\indent
76801 Saint-Etienne du Rouvray, France
}
\email{xuchaojiang@nuaa.edu.cn}

\date{\today}


\maketitle
\begin{abstract}
 In this work, we study the linear Landau equation with soft potential and show that the solution to the Cauchy problem with initial datum in $L^{2}(\mathbb{R}^3)$ enjoys an analytic regularizing effect, and the evolution of analytic radius is same as heat equations.
\end{abstract}

\section{Introduction}
In this work, we study the Cauchy problem of spatially homogeneous Landau equation :
\begin{equation}\label{eq1.10}
\left\{
\begin{array}{ll}
   \partial_t F= Q(F,F),\\
  F|_{t=0}=F_0,
\end{array}
\right.
\end{equation}
where $F=F(t,v)\ge0$ is the density distribution function depending on the velocity variables
$v\in\mathbb{R}^{3}$ and the time $t\geq0$. The Landau bilinear collision operator is given by
\begin{equation*}
Q(G,F)(v)
=\sum^3_{j, k=1}\partial_j\left(\int_{\mathbb{R}^{3}}
a_{jk}(v-v_*)
\big[G(v_*)(\partial_kF)(v)-(\partial_kG)(v_*)F(v)\big]
dv_*\right),
\end{equation*}
where
\begin{equation}\label{aij}
a_{jk}(v)=(\delta_{jk}|v|^2-v_jv_k)|v|^{\gamma},\quad\,\gamma\geq-3.
\end{equation}
One calls hard potentials if $\gamma>0$, Maxwellian molecules if $\gamma=0$, soft potentials if $\gamma\in ]-3,0[$ and Coulombian potential if $\gamma=-3$,
We shall study the linearization of the Landau equation \eqref{eq1.10} near the absolute Maxwellian distribution
$$
\mu(v)=(2\pi)^{-\frac 32}e^{-\frac{|v|^{2}}{2}}.
$$
Considering the fluctuation of density distribution function
$$
F(t,v)=\mu(v)+\sqrt{\mu}(v)f(t,v),
$$
since $Q(\mu,\mu)=0$, the Cauchy problem \eqref{eq1.10} is reduced to the Cauchy problem
\begin{equation} \label{eq-1}
\left\{ \begin{aligned}
         &\partial_t f+\mathcal{L}(f)+\mathcal{K}(f)={\bf \Gamma}(f, f),\,\,\, t>0,\, v\in\mathbb{R}^3,\\
         &f|_{t=0}=f_{0},
\end{aligned} \right.
\end{equation}
with $F_0(v)=\mu+\sqrt{\mu}f_0(v)$, where
\begin{align*}
{\bf \Gamma}(f, f)=\mu^{-\frac 12}Q_(\sqrt{\mu}f,\sqrt{\mu}f),
\end{align*}
and
\begin{equation}\label{l1}
\mathcal{L}f=-{\bf \Gamma}(\sqrt{\mu}, f),\qquad \mathcal{K}f=-{\bf \Gamma}(f,\sqrt{\mu}).
\end{equation}
In this work, we study Cauchy problem of the linear Landau equation,
such as
\begin{equation} \label{eq-1}
\left\{
\begin{array}{ll}
\partial_t f=\mathcal{L} f,\\
f|_{t=0}=f_0\in L^{2}(\mathbb{R}^{3}).
\end{array}
\right.
\end{equation}
Using the references \cite{Guo-2002} and \cite{WW}, we show that the diffusion part $\mathcal{L}$ is written as follows
\begin{equation}\label{l1-2}
\mathcal{L}f=-\nabla_v\cdot[A(v)\nabla_vf]+(A(v)\frac{v}{2}\cdot\frac{v}{2})f
-\nabla_v\cdot[A(v)\frac{v}{2}]f,
\end{equation}
with $A(v)=(\bar{a}_{ij})_{1\leq i,j\leq3}$ is a symmetric matrix where
\begin{align}\label{A}
\bar{a}_{ij}=a_{ij}*\mu=\int_{\mathbb{R}^3}(\delta_{jk}|v-v'|^2-(v_j-v'_j)(v_k-v'_k))
|v-v'|^{\gamma}\mu(v')d v'.
\end{align}

For the hard potential case, the existence, uniqueness of the solution to Cauchy problem for the spatially homogeneous Landau equation has already been treated in \cite{DV},\cite{Villani1998-2} under rather weak assumption on the initial datum. Moreover, they prove the smoothness of the solution in
$C^\infty(]0, +\infty[; \mathcal{S}(\mathbb{R}^3))$.  In \cite{ChenLiXu6}, Chen-Li-Xu improve this smoothing property and prove that the solution is in fact analytic for any $t>0$ (See  \cite{ChenLiXu3,ChenLiXu5} for the Gevrey regularity).

In the Maxwellian molecules case, in \cite{NYKC2}, Lerner,  Morimoto, Pravda-Starov and Xu study the spatially homogeneous non-cutoff Boltzmann equation and Landau equation in a close-to-equilibrium framework and show that the solution enjoys the Gelfand-Shilov smoothing effect (see also \cite{MPX}, \cite{LX-3} and \cite{Li2}). This implies that the nonlinear spatial homogeneous Landau equation has the same smoothing effect properties as the classic heat equation or harmonic oscillators heat equation.  In addition,  starting from a $L^2$ initial datum at $t=0$, the solution of Cauchy problem is spatial analytic for any $t>0$ and the analytic radius is $c_0 t^{\frac 12}$.  In the non-Maxwellian case, we can't use the Fourier transformation and spectral decomposition as in \cite{NYKC2,Li2,LX-3,MPX}.
Recently, Li and Xu in \cite{LX2022} proved the analytic smoothing effect of the solution to the nonlinear Landau equation with hard potentials and Maxwellian molecules case, that is $\gamma\geq0$.

In this work, we study the Cauchy problem of the linear Landau equation with soft potential $-3<\gamma<0$, we show the solution enjoys an analytic smoothing effect with the analytic radius $c_0 t^{\frac 12}$. The main theorem is in the following.

\begin{theorem}\label{trick}
For soft potential $-3< \gamma<0$ and for any $T>0$, the initial datum $f_0\in L^2(\mathbb{R}^3)$,  the Cauchy problem \eqref{eq-1} admits a unique weak solution
$$
f\in L^{\infty}([0,T]; L^2(\mathbb{R}^3)).
$$
Moreover, for any $\alpha\in \mathbb{N}^3$, $\tilde{t}=\min\{t, 1\}$, we have,
$$
\|\tilde{t}^{\frac{|\alpha|}{2}}\langle v\rangle^{\frac{\gamma|\alpha|}{2}}\partial^{\alpha}f(t)\|_{L^{2}(\mathbb{R}^3)}\leq C^{|\alpha|+1}\alpha!,
\qquad t\in [0, T].
$$
\end{theorem}

\begin{remark}\label{remark1}
Equivalently, for any $m\in \mathbb{N}$, we have,  for soft potential $-3<\gamma<0$,
\begin{align*}
&\|\tilde{t}^{\frac{m}{2}}\langle v\rangle^{\frac{\gamma m}{2}}\nabla^m f(t)\|^2_{L^2(\mathbb{R}^{3})}\\
&=\sum_{|\alpha|=m}\frac{(m!)^2}{(\alpha!)^2}\|t^{\frac{|\alpha|}{2}}
\langle v\rangle^{\frac{\gamma |\alpha|}{2}}\partial^{\alpha}f(t)\|^2_{L^2(\mathbb{R}^{3})}
\leq (3C)^{2m+2}(m!)^2,\qquad t\in [0, T].
\end{align*}
\end{remark}

This paper is arranged as follows : We prove the ultra-analytic for the coefficient of the Landau operator in Section \ref{S2}.  In Section \ref{S3}, we deal with estimation of the commutators and prove the coercivity property of the linear Landau operator.
In Section \ref{S4}, we study the Cauchy problem for linear Landau equation, and show that the existence and  uniqueness properties of the weak solution. The analytic smoothing effect of the weak solution for the linear Landau equation with soft potential will be proved in Section \ref{S5}.  In the Appendix \ref{Appendix}, we introduce the Hermite operator and related results.

\section{ultra-analytic for the coefficient of Landau operator}\label{S2}
For $\gamma\in \mathbb{R}$, denote
$$
\|\langle v\rangle^{\gamma}f\|_{L^{2}(\mathbb{R}^{3})}=\| f\|_{2,\gamma},\ \ \|\langle v\rangle^{\gamma}f\|_{L^{\infty}(\mathbb{R}^{3})}=\| f\|_{\infty,\gamma}
$$
where we use the notations $\langle v\rangle=(1+| v|^2)^{1/2}$.

In addition, for the matrix $A$ defined in \eqref{A}, we denote
\begin{equation}\label{{L2A}}
\| u\|_{A}^{2}=\sum_{j,k=1}^{3}\int _{\mathbb{R}^3}(\bar{a}_{jk}\partial_{j}u\partial_{k}u+\frac{1}{4}\bar{a}_{jk}v_{j}v_{k}u^{2})dv,
\end{equation}
and weighted norm, for $\theta\in \mathbb{R}$,
\begin{equation}\label{{L2AW}}
\| u\|_{A,\theta}^{2}=\sum_{j,k=1}^{3}\int _{\mathbb{R}^3}\langle v\rangle^{2\theta}(\bar{a}_{jk}\partial_{j}u\partial_{k}u+\frac{1}{4}\bar{a}_{jk}v_{j}v_{k}u^{2})dv.
\end{equation}

From formula $(21)$ of Corollary 1 in \cite{Guo-2002}, for any $\theta\in \mathbb{R}$, there exist $C_{1}>0$, such that
\begin{equation}\label{definition}
\| f\|_{A,\theta}^{2}\geq C_{1}(\|P_{v}\triangledown f\|_{2,\frac{\gamma}{2}+\theta}^{2}+\|(I-P_{v})\triangledown f\parallel_{2,1+\frac{\gamma}{2}+\theta}^{2}+\| f\|_{2,1+\frac{\gamma}{2}+\theta}^{2}).
\end{equation}
where for any vector-valued function $G(v)=(G_{1},G_{2},G_{3})$, we define the projection to the vector $v=(v_{1},v_{2},v_{3})$ as
$$
P_{v}G_{i}=\sum_{j=1}^{3}G_{j}v_{j}\frac{v_{i}}{| v|^{2}},\ \  1\leq i\leq3.
$$
Notice that $\nabla f=P_{v}\nabla f+(I-P_{v})\nabla f$, we have
\begin{equation}\label{definition2}
\|f\|_{A,\theta}\geq C_{1}(\|\nabla f\|_{2,\frac{\gamma}{2}+\theta}+\|f\|_{2,1+\frac{\gamma}{2}+\theta}).
\end{equation}
We can also refer to \cite{DeL} and their references. Remark that
the weighted of $f$ and $\nabla f$ are different in the norm $L_{A}^{2}$,

Firstly, for any $\gamma>-3$ and $\delta>0$, we have
\begin{equation}\label{cov}
\int_{\mathbb{R}^3}|v-w|^{\gamma}e^{-\delta|w|^2}dw\lesssim \langle v\rangle^{\gamma}.
\end{equation}
Which imply
\begin{equation}\label{cov-2}
|\bar{a}_{i j}(v)|\lesssim \langle v\rangle^{\gamma+2}.
\end{equation}

In the following, we prove that the coefficients of linear Landau operator are ultra-analytic.

\begin{lemma}\label{ultra-analytic}
For any $\beta\in \mathbb{N}^3$ with $|\beta|\geq1$ and  $\bar{a}_{ij}$ was defined in \eqref{A} with $-3<\gamma<0$ , then we have,
\begin{equation}\label{ultra-1}
|\partial^{\beta}\bar{a}_{ij}(v)|
\lesssim\langle v\rangle^{\gamma+1}\sqrt{\beta!}.
\end{equation}
Moreover, for any $\beta\in \mathbb{N}^3$,
\begin{equation}\label{ultra-2}
\begin{split}
&|\partial^{\beta}(\sum_{i,j=1}^{3}
\partial_ia_{ij}\ast(v_{j}\mu))|
\lesssim \langle v\rangle^{\gamma+1}(|\beta|+1)\sqrt{\beta!},\\
&|\partial^{\beta}(\sum_{i,j=1}^{3}
\bar{a}_{ij}v_{i}v_{j})|\lesssim \langle v\rangle^{\gamma+1}(|\beta|+1)\sqrt{\beta!}.
\end{split}
\end{equation}
\end{lemma}
\begin{proof}
For $\beta\in \mathbb{N}^3$ with $|\beta|\geq1$, without loss of generality, we set
$$\beta_1=\max(\beta_1,\beta_2,\beta_3).$$
  Notice that $\bar{a}_{ij}=a_{ij}\ast\mu$, then
\begin{align*}
\partial^{\beta}\bar{a}_{ij}=\partial^{\beta}[a_{ij}\ast\mu]
=\partial_{1}a_{ij}\ast\partial^{\beta-e_1}\mu
\end{align*}
Direct calculation shows that, for any $1\leq i\leq3$,
\begin{equation}\label{P-A}
\partial_i(\sqrt{\mu}f)=\sqrt{\mu}(\partial_i-\frac{v_i}{2})f=-\sqrt{\mu}A_{+,i}f.
\end{equation}
For more details of the operators $A_{\pm,i}$, we can refer to the Appendix \ref{Appendix}.  By using the fact
$$|\partial_1a_{ij}(v)|\lesssim |v|^{\gamma+1},$$
it follows from Cauchy-Schwartz's inequality and \eqref{cov} that
\begin{align*}
|\partial^{\beta}\bar{a}_{ij}(v)|
&\leq\int_{\mathbb{R}^3}|\partial_1a_{ij}(v-w)(\sqrt{\mu}
A_{+,1}^{\beta_{1}-1}A_{+,2}^{\beta_{2}}A_{+,3}^{\beta_{3}}\sqrt{\mu})(w)|dw\\
&\lesssim\left(\int_{\mathbb{R}^3}|v-w|^{\frac{3}{2}(\gamma+1)}|\sqrt{\mu}(w)dw\right)^{\frac{2}{3}}
\|A_{+,1}^{\beta_{1}-1}A_{+,2}^{\beta_{2}}
A_{+,3}^{\beta_{3}}\sqrt{\mu}\|_{L^3(\mathbb{R}^3)}\\
&\lesssim(1+| v|)^{\gamma+1}\|A_{+,1}^{\beta_{1}-1}A_{+,2}^{\beta_{2}}
A_{+,3}^{\beta_{3}}\sqrt{\mu}\|_{L^3(\mathbb{R}^3)}
\end{align*}
where we use the fact
$$\frac{3(\gamma+1)}{2}>-3,\quad
 \forall
 -3<\gamma<0.$$
 Due to $\sqrt{\mu}(v)=\Psi_{0}(v)$, we can deduce from \eqref{H4} that
\begin{equation}\label{psi}
A_{+,1}^{\beta_{1}-1}A_{+,2}^{\beta_{2}}A_{+,3}^{\beta_{3}}\sqrt{\mu}
=\sqrt{(\beta_{1}-1)!\beta_{2}!\beta_{3}!}\Psi_{\beta-e_1},
\end{equation}
where $\{\Psi_{\alpha}\}_{\alpha\in \mathbb{N}^3}$ is the orthonormal basis in $L^2(\mathbb{R}^3)$.
The H\"older's inequality and Poincar\'e's inequality implies
$$\|f\|^2_{L^3(\mathbb{R}^3)}\leq \|f\|_{L^2(\mathbb{R}^3)}\|f\|_{L^6(\mathbb{R}^3)}\lesssim \|f\|_{L^2(\mathbb{R}^3)}\|\nabla f\|_{L^2(\mathbb{R}^3)},$$
along with the equalities \eqref{H3} and \eqref{H4} shows that,
\begin{align*}
|\partial^{\beta}\bar{a}_{ij}(v)|
&\lesssim \langle v\rangle^{\gamma+1}\|\nabla\Psi_{\beta-e_1}\|^{\frac{1}{2}}_{L^2(\mathbb{R}^3)}\\
&\lesssim \langle v\rangle^{\gamma+1}\sum^3_{k=1}\sqrt{(\beta-e_1+e_k)!}
\end{align*}
Consider that $\beta_1=\max(\beta_1,\beta_2,\beta_3)$,
we get that
$$
|\partial^{\beta}\bar{a}_{ij}(v)|
\lesssim\langle v\rangle^{\gamma+1}\sqrt{\beta!}.
$$
For the estimate of the remaining inequalities \eqref{ultra-2}, an integration by parts inside the convolution show that,
\begin{align*}
\sum_{i,j=1}^{3}\partial_ia_{ij}\ast(v_{j}\mu)
&=-\sum_{i,j=1}^{3}\partial_ia_{ij}\ast(\sqrt{\mu}A_{+,j}\sqrt{\mu}),\\
\sum_{i,j=1}^{3}a_{ij}\ast(v_iv_{j}\mu)
&=\sum_{i,j=1}^{3}a_{ij}\ast(\sqrt{\mu}(A_{+,i}A_{+,j}\sqrt{\mu}
-\delta_{i,j}A_{+,j}\sqrt{\mu}))\\
&=-\sum_{i,j=1}^{3}\partial_ja_{ij}\ast(\sqrt{\mu}(A_{+,i}\sqrt{\mu}
-\delta_{i,j}\sqrt{\mu})).
\end{align*}
By using the fact
$$|\partial_ia_{ij}(v)|\lesssim |v|^{\gamma+1},\quad|\partial_ja_{ij}(v)|\lesssim |v|^{\gamma+1}.$$
A calculation similar to that as above, it follows that
\begin{align*}
&|\partial^{\beta}(\sum_{i,j=1}^{3}
\partial_ia_{ij}\ast(v_{j}\mu))|\\
&\lesssim\sum_{j=1}^{3}\left(\int_{\mathbb{R}^3}
|v-w|^{\frac{3}{2}(\gamma+1)}e^{-\frac{|w|^2}{2}}dw\right)^{\frac{1}{2}}\sqrt{(\beta+e_j)!}\|\Psi_{\beta+e_j}\|_{L^3(\mathbb{R}^3)}\\
&\lesssim \langle v\rangle^{\gamma+1}
\sum^3_{j=1}\sqrt{(\beta+e_j)!}\|\Psi_{\beta+e_j}\|_{L^3(\mathbb{R}^3)}\\
&\lesssim \langle v\rangle^{\gamma+1}\sum^3_{j=1}\sqrt{(\beta+e_j)!}\|\nabla\Psi_{\beta+e_j}\|^{\frac{1}{2}}_{L^2(\mathbb{R}^3)}\\
&\lesssim \langle v\rangle^{\gamma+1}\sum^3_{k,j=1}\sqrt{(\beta+e_j+e_k)!}\\
&\leq \langle v\rangle^{\gamma+1}(|\beta|+1)\sqrt{\beta!}.
\end{align*}
The same estimate holds true for the last term, such that
\begin{align*}
|\partial^{\beta}(\sum_{i,j=1}^{3}
\bar{a}_{ij}v_iv_{j})|
&\lesssim \langle v\rangle^{\gamma+1}(|\beta|+1)\sqrt{\beta!}.
\end{align*}
We end the proof of Lemma \ref{ultra-analytic}.
\end{proof}

 In order to prove the coercivity of the linear Landau operator, we need one more estimate to control the weighted $\bar{a}_{ij}$.
\begin{lemma}\label{1-estimate}
For $f,g\in\mathcal{S}(\mathbb{R}^3)$,  for any  $\beta\in \mathbb{N}^3$ and $\theta\in \mathbb{R}$, we have
\begin{equation}\label{estimate1}
|\sum_{i,j=1}^{3}\langle\langle v\rangle^{2\theta}\partial^{\beta}\bar{a}_{ij}\partial_{j}f,
\partial_{i}g\rangle|\lesssim \sqrt{\beta!}\| f\|_{A,\theta}\| g\|_{A,\theta}.
\end{equation}
where $\bar{a}_{ij}$ was defined in \eqref{A} with $-3< \gamma<0$.
\end{lemma}
\begin{proof}
In fact, the inner product
\begin{equation*}
\langle\langle v\rangle^{2\theta}\partial^{\beta}\bar{a}_{ij}\partial_{j}f,\partial_{i}g\rangle=\iint_{\mathbb{R}^3\times \mathbb{R}^3}\langle v\rangle^{2\theta}a_{ij}(v-w)\partial^{\beta}\mu(w)\partial_{j}f(v)\partial_{i}g(v)dwdv.
\end{equation*}
We decompose the integration region $[v,w]\in \mathbb{R}^3\times \mathbb{R}^3$ into three parts:
$$
\{| v|\leq1\},\ \  \{2| w|\geq| v|,| v|\geq1\},\ \ and\ \  \{2| w|\leq| v|,| v|\geq1\}.
$$
For the first part $\{| v|\leq1\}$, by Lemma \ref{ultra-analytic}, we have
\begin{align*}
&\left|\sum_{i,j=1}^{3}\iint_{\{| v|\leq1\}\times\mathbb{R}^3}\langle v\rangle^{2\theta}a_{ij}(v-w)\partial^{\beta}\mu(w)\partial_{j}f(v)\partial_{i}g(v)dwdv\right|\\
&\lesssim\sum_{i,j=1}^{3}\int_{| v|\leq1}\langle v\rangle^{\gamma+2}\sqrt{\beta!}
\langle v\rangle^{2\theta}|\partial_{j}f\partial_{i}g| dv\\
&\leq\sqrt{\beta!}\|\langle v\rangle^{\frac{\gamma}{2}}\nabla f\|_{2,\theta}\|\langle v\rangle^{\frac{\gamma}{2}}\nabla g\|_{2,\theta}\lesssim \sqrt{\beta!}\|f\|_{A,\theta}\|g\|_{A,\theta}.
\end{align*}
For the second part $\{2| w|\geq| v|,| v|\geq1\}$, we have
$$e^{-\frac{|w|^2}{4}}\leq e^{-\frac{|w|^2}{8}} e^{-\frac{|v|^2}{32}}.$$
Similar to the proof as Lemma \ref{ultra-analytic}, one can verify that
\begin{align*}
&\left|\sum_{i,j=1}^{3}\iint_{\{2| w|\geq| v|,| v|\geq1\}}\langle v\rangle^{2\theta}a_{ij}(v-w)\sqrt{\mu(w)}\Psi_{\beta}(w)\partial_{j}f(v)\partial_{i}g(v)dwdv\right|\\
&\lesssim\sum_{i,j=1}^{3}\int_{\mathbb{R}^3}\langle v\rangle^{\gamma+2}e^{-\frac{|v|^2}{32}}\sqrt{\beta!}
\langle v\rangle^{2\theta}|\partial_{j}f\partial_{i}g| dv\\
&\lesssim\sqrt{\beta!}\|f\|_{A,\theta}\|g\parallel_{A,\theta}.
\end{align*}
We finally consider the third part $\{2| w|\leq| v|,| v|\geq1\}.$ Expanding $a_{ij}(v-w)$ to get
\begin{equation*}
a_{ij}(v-w)=a_{ij}(v)+\sum^3_{k=1}\partial_{k}a_{ij}(v)w_k+\frac{1}{2}\sum_{k,l=1}^{3}\left(\int_{0}^{1}\partial_{kl}a_{ij}(v-sw)ds\right)w_{k}w_l.
\end{equation*}
Along with the fact that
$$
\sum_{i=1}^{3}a_{ij}v_{i}=0, \quad \sum_{i,j=1}^{3}\partial_{k}a_{ij}(v)v_iv_j=-2\sum_{j=1}^{3}a_{kj}(v)v_j=0,
$$
show immediately
\begin{align*}
&\sum_{i,j=1}^{3}\iint_{\{2| w|\leq| v|,| v|\geq1\}}\langle v\rangle^{2\theta}a_{ij}(v-w)\partial^{\beta}\mu(w)\partial_{j}f(v)\partial_{i}g(v)dwdv\\
=&\sum_{i,j=1}^{3}\iint_{2| w|\leq| v|,| v|\geq1}\langle v\rangle^{2\theta}\partial^{\beta}\mu(w)a_{ij}(v)[(I-P_{v})\partial_{j}f][(I-P_{v})\partial_{i}g]dwdv\\
&+\sum_{i,j=1}^{3}\sum^3_{k=1}\iint_{2| w|\leq| v|,| v|\geq1}\langle v\rangle^{2\theta}w_k\partial^{\beta}\mu(w)\partial_{k}a_{ij}(v)[P_{v}\partial_{j}f][(I-P_{v})\partial_{i}g]dwdv\\
&+\sum_{i,j=1}^{3}\sum^3_{k=1}\iint_{2| w|\leq| v|,| v|\geq1}\langle v\rangle^{2\theta}w_k\partial^{\beta}\mu(w)\partial_{k}a_{ij}(v)[(I-P_{v})\partial_{j}f][\partial_{i}g]dwdv\\
&+\frac{1}{2}\sum_{k,l=1}^{3}\sum_{i,j=1}^{3}\int_{0}^{1}\iint_{2| w|\leq| v|,| v|\geq1}\langle v\rangle^{2\theta}w_{k}w_l\partial^{\beta}\mu(w)\partial_{kl}a_{ij}(v-sw)\partial_{j}f\partial_{i}gdwdvds.
\end{align*}
Since $2| w|\leq| v|,| v|\geq1,0<s<1,$ for $-3<\gamma<0$, we have
\begin{align*}
&|a_{ij}(v)|\lesssim |v|^{\gamma+2}\lesssim\langle v\rangle^{\gamma+2};\quad |\partial_{k}a_{ij}(v)|\lesssim \langle v\rangle^{\gamma+1};\\
&|\partial_{kl}a_{ij}(v-sw)|\leq C| v-sw|^{\gamma}\leq C4^{-\gamma}(1+| v|)^{\gamma}.
\end{align*}
It follows from the inequality \eqref{psi} and the norm equality \eqref{definition}  that
\begin{align*}
&\left|\sum_{i,j=1}^{3}\iint_{\{2| w|\leq| v|,| v|\geq1\}}\langle v\rangle^{2\theta}a_{ij}(v-w)\partial^{\beta}\mu(w)\partial_{j}f(v)\partial_{i}g(v)dwdv\right|\\
&\lesssim\sqrt{\beta!}\|(I-P_{v})\nabla f\|_{2,\frac{\gamma+2}{2}+\theta}
\|(I-P_{v})\nabla g\|_{2,\frac{\gamma+2}{2}+\theta}\\
&\quad+\sqrt{\beta!}
\|(I-P_{v})\nabla f\|_{2,\frac{\gamma+2}{2}+\theta}
\|P_{v}\nabla g\|_{2,\frac{\gamma}{2}+\theta}\\
&\quad+\sqrt{\beta!}\|(I-P_{v})\nabla f\|_{2,\frac{\gamma+2}{2}+\theta}
\|\nabla g\|_{2,\frac{\gamma}{2}+\theta}+\sqrt{\beta!}\|\nabla f\|_{2,\frac{\gamma}{2}+\theta}
\|\nabla g\|_{2,\frac{\gamma}{2}+\theta}\\
&\lesssim\sqrt{\beta!}\| f\|_{A,\theta}\|g\|_{A,\theta}.
\end{align*}
This is the inequality \eqref{estimate1}.
  We end the proof of Lemma \ref{1-estimate}.
\end{proof}

\section{Estimations of commutators }\label{S3}

\begin{proposition}\label{multiply-estimate-weighted}
Let $f\in\mathcal{S}(\mathbb{R}^3)$,  $\mathcal{L}$ was defined in \eqref{eq-1}, for any $\alpha\in\mathbb{N}^3$ and $\theta\in \mathbb{R}$, there exist a positive constant $C_{0}>0$ which is independent on $\alpha$ and $\theta$, such that,
\begin{align*}
\Big(\langle v\rangle^{2\theta}\partial^{\alpha}\mathcal{L} f,& \partial^{\alpha}f\Big)_{L^{2}(\mathbb{R}^{3})}
\leq-\| \partial^{\alpha}f\|_{A,\theta}^{2}+C_0
\langle\theta\rangle\|\partial^{\alpha} f\|_{A,\theta}\|\partial^{\alpha}f \|_{2,\frac{\gamma}{2}+\theta}\\
&+C_0\sum_{\substack{|\beta|\ge 1\\
\beta\leq\alpha}}C_{\alpha}^{\beta} \sqrt{\beta!} \|\partial^{\alpha-\beta}f\|_{A,\theta}
\left(\|\partial^{\alpha}f\|_{A,\theta}+|\theta|\|\partial^{\alpha}f\|_{2,\theta+\frac{\gamma}{2}}\right)\\
&+C_0\sum_{\substack{|\beta|\geq1\\
\beta\leq\alpha}}C_{\alpha}^{\beta}|\beta|\sqrt{\beta!}
\|\partial^{\alpha-\beta} f\|_{A,\theta}\|\partial^{\alpha}f\|_{2,\theta+\frac{\gamma}{2}}.
\end{align*}
\end{proposition}
\begin{remark} Remark that, we have
\begin{enumerate}
  \item For $\alpha=0, \theta=0$, we have
  \begin{equation}\label{Remark-1}
\Big(\mathcal{L} f, f\Big)_{L^{2}(\mathbb{R}^{3})}
\leq-\frac 12\|f\|_{A}^{2}+C_0\|f\|^2_{\frac{\gamma}{2}}.
\end{equation}
  \item For $|\alpha|\ge 1$,
\begin{equation}\label{Remark-2}
\|\partial^{\alpha}f\|_{2,\frac{\gamma}{2}+\theta}\lesssim \|\nabla\partial^{\alpha-e_{j_0}} f\|_{2,\frac{\gamma}{2}+\theta}\lesssim \|\partial^{\alpha-e_{j_0}}f\|_{A,\theta},
\end{equation}
where $\alpha_{j_0}=\max\{\alpha_1, \alpha_2, \alpha_3\}$.
\end{enumerate}
\end{remark}

\begin{proof}
Recalled the formula $\mathcal{L}f$ in \eqref{eq-1}, for the smooth function $f$, integrated by parts, we have
\begin{align*}
-\left(\langle v\rangle^{2\theta}\partial^{\alpha}\mathcal{L}f,\partial^{\alpha}f\right)_{L^{2}(\mathbb{R}^{3})}
=&\sum_{i,j=1}^{3}\int_{\mathbb{R}^3}\langle v\rangle^{2\theta}\left\{\partial^{\alpha}
(\bar{a}_{ij}\partial_{j}f)\right\}
\partial^{\alpha}\partial_{i}f dv\\
&+\sum_{i,j=1}^{3}\int_{\mathbb{R}^3}\left\{\partial_i\langle v\rangle^{2\theta}\right\}\left\{\partial^{\alpha}
(\bar{a}_{ij}\partial_{j}f)\right\}
\partial^{\alpha}f dv\\
&+\frac{1}{4}\sum_{i,j=1}^{3}
\int_{\mathbb{R}^3}\langle v\rangle^{2\theta}\left\{\partial^{\alpha}(\bar{a}_{ij}v_{i}v_{j}f)\right\}
\partial^{\alpha}f dv\\
&-\frac{1}{2}\sum_{i,j=1}^{3}\int_{\mathbb{R}^3}\langle v\rangle^{2\theta}
\left\{\partial^{\alpha}\partial_{i}((\bar{a}_{ij}v_{j})f)\right\}\partial^{\alpha}fdv.
\end{align*}
Then, using \eqref{{L2AW}}, we have,
\begin{equation}\label{linear-decomp}
-\left(\langle v\rangle^{2\theta}\partial^{\alpha}\mathcal{L}f,\partial^{\alpha}f\right)_{L^{2}(\mathbb{R}^{3})}
=\| \partial^{\alpha}f\|_{A,\theta}^{2}+\mathbf{R}(f),
\end{equation}
with
\begin{align*}
\mathbf{R}(f)
=&\sum_{i,j=1}^{3}\int_{\mathbb{R}^3}\langle v\rangle^{2\theta}\left\{[\partial^{\alpha}
, \bar{a}_{ij}]\partial_{j}f)\right\}
\partial^{\alpha}\partial_{i}f dv\\
&+\sum_{i,j=1}^{3}\int_{\mathbb{R}^3}\left\{\partial_i\langle v\rangle^{2\theta}\right\}\left\{\partial^{\alpha}
(\bar{a}_{ij}\partial_{j}f)\right\}
\partial^{\alpha}f dv\\
&+\frac{1}{4}\sum_{i,j=1}^{3}
\int_{\mathbb{R}^3}\langle v\rangle^{2\theta}\left\{[\partial^{\alpha}, \bar{a}_{ij}v_{i}v_{j}] f)\right\}
\partial^{\alpha}f dv\\
&-\frac{1}{2}\sum_{i,j=1}^{3}\int_{\mathbb{R}^3}\langle v\rangle^{2\theta}
\left\{\partial^{\alpha}\partial_{i}((\bar{a}_{ij}v_{j})f)\right\}\partial^{\alpha}fdv.
\end{align*}
By using Leibniz formula,
$$
\partial^{\alpha}(g h)=g \partial^{\alpha}h+\sum_{\substack{|\beta|=1\\
\beta\leq\alpha}}C_{\alpha}^{\beta}\partial^{\beta}g
\partial^{\alpha-\beta}h+\sum_{\substack{|\beta|\geq2\\
\beta\leq\alpha}}C_{\alpha}^{\beta}\partial^{\beta}g
\partial^{\alpha-\beta}h,
$$
and
$$
\bar{a}_{ij}=a_{ij}\ast \mu,\quad
\sum_{j}\bar{a}_{ij}v_{j}=\sum_{j}a_{ij}\ast (v_{j}\mu),
\quad
\sum_{i,j}\bar{a}_{ij}v_i v_{j}=\sum_{j}a_{ij}\ast (v_iv_{j}\mu),
$$
it follows that
$$
\mathbf{R}(f)
=\mathbf{R}_{0}(f)+\mathbf{R}_{1}(f)
$$
where
\begin{align*}
\mathbf{R}_{0}(f)=&\sum_{i,j=1}^{3}\int_{\mathbb{R}^3}(\partial_i\langle v\rangle^{2\theta})
(\bar{a}_{ij}\partial_{j}\partial^{\alpha}f)
\partial^{\alpha}f dv\\
&-\frac{1}{2}\sum_{i,j=1}^{3}
\int_{\mathbb{R}^3}\langle v\rangle^{2\theta}(\partial_{i}a_{ij})\ast(v_{j}\mu)|\partial^{\alpha}f|^{2}dv
=\mathbf{R}_{01}(f)+\mathbf{R}_{02}(f);
\end{align*}
\begin{align*}
\mathbf{R}_{1}(f)=&\sum_{\substack{|\beta|\ge1\\
\beta\leq\alpha}}C_{\alpha}^{\beta} \sum_{i,j=1}^{3}\int_{\mathbb{R}^3}\langle v\rangle^{2\theta}
(\partial^\beta\bar{a}_{ij})(\partial^{\alpha-\beta}\partial_{j}f)
\partial^{\alpha}\partial_{i}fdv\\
&+\sum_{\substack{|\beta|\ge1\\
\beta\leq\alpha}}C_{\alpha}^{\beta}  \sum_{i,j=1}^{3}\int_{\mathbb{R}^3}\left\{\partial_i\langle v\rangle^{2\theta}\right\}
(\partial^\beta\bar{a}_{ij})\left\{\partial_{j}\partial^{\alpha-\beta}f\right\}
\partial^{\alpha}fdv
\\
&+\frac{1}{4}\sum_{\substack{|\beta|\ge1\\
\beta\leq\alpha}}C_{\alpha}^{\beta}  \sum_{i,j=1}^{3}\int_{\mathbb{R}^3}\langle v\rangle^{2\theta}
\left\{\partial^\beta(a_{ij}\ast(v_{i}v_{j}\mu))\right\}\partial^{\alpha-\beta}f
\partial^{\alpha}f dv\\
&-\frac{1}{2}\sum_{\substack{|\beta|\ge1\\
\beta\leq\alpha}}C_{\alpha}^{\beta}  \sum_{i,j=1}^{3}\int_{\mathbb{R}^3}\langle v\rangle^{2\theta}(\partial_{i}\partial^\beta a_{ij})\ast(v_{j}\mu)\partial^{\alpha-\beta}f\partial^{\alpha}fdv\\
=&\mathbf{R}_{11}+\mathbf{R}_{12}+\mathbf{R}_{13}+\mathbf{R}_{14}.
\end{align*}
So that the proof of the Proposition \ref{multiply-estimate-weighted} is reduced to the estimations of terms $\mathbf{R}_{0}(f)$ and $\mathbf{R}_{1}(f)$ , which will give by the following two lemmas.
\end{proof}

\begin{lemma}\label{lemma3.1}
We have, for $\alpha\in\mathbb{N}^3$,
\begin{equation}\label{R-0}
|\mathbf{R}_{0}(f)|\lesssim {\color{blue}\langle\theta\rangle}|\|\partial^{\alpha} f\|_{A,\theta}\|\partial^{\alpha}f \|_{2,\frac{\gamma}{2}+\theta}.
\end{equation}
\end{lemma}
\begin{proof}
By using \eqref{cov}, we get
$$
|(\partial_{i}a_{ij})\ast(v_{j}\mu)| \lesssim \langle v\rangle^{\gamma+1},\quad
|(a_{ij})\ast(v_{j}\mu)| \lesssim \langle v\rangle^{\gamma+2}
$$
then \eqref{definition2} imply
\begin{align*}
|\mathbf{R}_{02}(f)|&=\left|\sum_{i,j=1}^{3}
\int_{\mathbb{R}^3}\langle v\rangle^{2\theta}(\partial_{i}a_{ij})\ast(v_{j}\mu)|\partial^{\alpha}f|^{2}dv\right|\\
&\lesssim\int_{\mathbb{R}^3}\langle v\rangle^{\gamma+1+2\theta}
|\partial^{\alpha}f|^{2}dv\\
&\lesssim\|\partial^{\alpha}f\|_{2,\frac{\gamma}{2}+\theta}
\|\partial^{\alpha}f\|_{2,\frac{\gamma}{2}+1+\theta} \lesssim\|\partial^{\alpha}f\|_{2,\frac{\gamma}{2}+\theta}\|\partial^{\alpha}f\|_{A,\theta}.
\end{align*}
For the term $\mathbf{R}_{01}$,  we use
$$
\partial_i\langle v\rangle^{2\theta}=2\theta\langle v\rangle^{2\theta-2}v_i,
$$
then,
\begin{align*}
\mathbf{R}_{01}(f)&=\sum_{i,j=1}^{3}\int_{\mathbb{R}^3}(\partial_i\langle v\rangle^{2\theta})
(\bar{a}_{ij}\partial_{j}\partial^{\alpha}f)
\partial^{\alpha}f dv\\
&=2\theta\sum_{i,j=1}^{3}\int_{\mathbb{R}^3}\langle v\rangle^{2\theta-2}
(\bar{a}_{ij}v_i\partial_{j}\partial^{\alpha}f)
\partial^{\alpha}f dv\\
&=2\theta\sum_{i,j=1}^{3}\int_{\mathbb{R}^3}\langle v\rangle^{2\theta-2}
({a}_{ij}\ast(v_{i}\mu)\partial_{j}\partial^{\alpha}f)
\partial^{\alpha}f dv.
\end{align*}
It follows that
\begin{align*}
|\mathbf{R}_{01}(f)|\leq&2|\theta|\sum_{j=1}^{3}\int_{\mathbb{R}^3}\langle v\rangle^{2\theta+\gamma}
|\partial_{j}\partial^{\alpha} f|
|\partial^{\alpha}f |dv\\
\lesssim
&|\theta|\sum_{j=1}^{3}\|\partial_{j}\partial^{\alpha}f\|_{2,\frac{\gamma}{2}+\theta}
\|\partial^{\alpha}f\|_{2,\frac{\gamma}{2}+\theta}\\
\lesssim& |\theta|\|\partial^{\alpha}f\|_{A,\theta}\|\partial^{\alpha}f\|_{2,\frac{\gamma}{2}+\theta}.
\end{align*}
Which give the estimate of $\mathbf{R}_{0}(f)$.
\end{proof}

\begin{lemma}\label{lemma3.2}
We have, for $\alpha\in\mathbb{N}^3, |\alpha|\ge 1$,
\begin{equation}\label{R-1}
|\mathbf{R}_{1}(f)|\lesssim \sum_{\substack{|\beta|=1\\
\beta\leq\alpha}}C_{\alpha}^{\beta}  \|\partial^{\alpha-\beta}f\|_{A,\theta}
\left(\|\partial^{\alpha}f\|_{A,\theta}+|\theta|\|\partial^{\alpha}f\|_{2,\theta+\frac{\gamma}{2}}\right).
\end{equation}
\end{lemma}

\begin{proof}
Now we estimate $\mathbf{R}_{1}$.  For the term $\mathbf{R}_{11}(f)$, by using the inequality \eqref{estimate1} in Lemma \ref{1-estimate} directly, we obtain that
\begin{align*}
|\mathbf{R}_{11}(f)|&=
\left|\sum_{\substack{|\beta|\ge1\\
\beta\leq\alpha}}C_{\alpha}^{\beta} \sum_{i,j=1}^{3}\int_{\mathbb{R}^3}\langle v\rangle^{2\theta}
(\partial^\beta\bar{a}_{ij})(\partial^{\alpha-\beta}\partial_{j}f)
\partial^{\alpha}\partial_{i}fdv\right|\\
&\lesssim
\sum_{\substack{|\beta|\ge1\\
\beta\leq\alpha}}C_{\alpha}^{\beta}\sqrt{\beta!}  \|\partial^{\alpha-\beta}f\|_{A,\theta}
\|\partial^{\alpha}f\|_{A,\theta}.
\end{align*}
For the two terms $\mathbf{R}_{13}(f), \mathbf{R}_{14}(f)$,  we can deduce from the inequality \eqref{ultra-2} in Lemma \ref{ultra-analytic} and \eqref{definition2} that
\begin{align*}
|\mathbf{R}_{13}(f)|+|\mathbf{R}_{14}(f)|&=\left|
\frac{1}{4}\sum_{\substack{|\beta|\ge1\\
\beta\leq\alpha}}C_{\alpha}^{\beta}  \sum_{i,j=1}^{3}\int_{\mathbb{R}^3}\langle v\rangle^{2\theta}
\left\{\partial^\beta(a_{ij}\ast(v_{i}v_{j}\mu))\right\}\partial^{\alpha-\beta}f
\partial^{\alpha}f dv\right|\\
&\quad +\left|\frac{1}{2}\sum_{\substack{|\beta|\ge1\\
\beta\leq\alpha}}C_{\alpha}^{\beta}  \sum_{i,j=1}^{3}\int_{\mathbb{R}^3}\langle v\rangle^{2\theta}(\partial_{i}\partial^\beta a_{ij})\ast(v_{j}\mu)\partial^{\alpha-\beta}f\partial^{\alpha}fdv
\right|\\
&\lesssim\sum_{\substack{|\beta|\ge1\\
\beta\leq\alpha}}C_{\alpha}^{\beta} \int_{\mathbb{R}^3}\langle v\rangle^{2\theta+\gamma+1}|\beta|\sqrt{\beta!}|\partial^{\alpha-\beta}f||\partial^{\alpha}f|dv\\
&\lesssim\sum_{\substack{|\beta|\ge1\\
\beta\leq\alpha}}C_{\alpha}^{\beta} |\beta|\sqrt{\beta!}\|\partial^{\alpha-\beta}f\|_{2,\theta+\frac{\gamma+2}{2}}
\|\partial^{\alpha}f\|_{2,\theta+\frac{\gamma}{2}}\\
&\lesssim\sum_{\substack{|\beta|\ge1\\
\beta\leq\alpha}}C_{\alpha}^{\beta} |\beta|\sqrt{\beta!}\|\partial^{\alpha-\beta}f\|_{A,\theta}
\|\partial^{\alpha}f\|_{2,\theta+\frac{\gamma}{2}}.
\end{align*}
For the term $\mathbf{R}_{12}(f)$, for $|\beta|\geq1,$ it follows from \eqref{ultra-1} in Lemma \ref{ultra-analytic} that
\begin{align*}
|\partial_i\langle v\rangle^{2\theta}
\partial^{\beta}\bar{a}_{ij}|
\lesssim |\theta|\langle v\rangle^{2\theta-1}\langle v\rangle^{\gamma+1}\sqrt{\beta!}\leq |\theta|\langle v\rangle^{2\theta+\gamma}\sqrt{\beta!} .
\end{align*}
Then we have
\begin{align*}
|\mathbf{R}_{12}(f)|&=
\left|\sum_{\substack{|\beta|\ge1\\
\beta\leq\alpha}}C_{\alpha}^{\beta}  \sum_{i,j=1}^{3}\int_{\mathbb{R}^3}\left\{\partial_i\langle v\rangle^{2\theta}\right\}
(\partial^\beta\bar{a}_{ij})\left\{\partial_{j}\partial^{\alpha-\beta}f\right\}
\partial^{\alpha}fdv\right|\\
&\lesssim|\theta|\sum_{\substack{|\beta|\ge1\\
\beta\leq\alpha}}C_{\alpha}^{\beta} \sqrt{\beta!}\sum^3_{j=1}\|\partial^{\alpha-\beta}\partial_jf\|_{2,\frac{\gamma}{2}+\theta}
\|\partial^{\alpha}f\|_{2,\frac{\gamma}{2}+\theta}\\
&\lesssim|\theta|\sum_{\substack{|\beta|\ge1\\
\beta\leq\alpha}}C_{\alpha}^{\beta}\sqrt{\beta!} \|\partial^{\alpha-\beta} f\|_{A,\theta}
\|\partial^{\alpha}f\|_{2,\frac{\gamma}{2}+\theta}.
\end{align*}
We conclude that
\begin{align*}
|\mathbf{R}_{1}(f)|\lesssim&
\sum_{\substack{|\beta|\geq1\\
\beta\leq\alpha}}C_{\alpha}^{\beta}\sqrt{\beta!}
\|\partial^{\alpha-\beta}f\|_{A,\theta}\|\partial^{\alpha}f\|_{A,\theta}\\
&+|\theta|\sum_{\substack{|\beta|\geq1\\
\beta\leq\alpha}}C_{\alpha}^{\beta}\sqrt{\beta!}
\|\partial^{\alpha-\beta}f\|_{A,\theta}\|\partial^{\alpha} f\|_{2,\theta+\frac{\gamma}{2}}\\
&+\sum_{\substack{|\beta|\geq1\\
\beta\leq\alpha}}C_{\alpha}^{\beta}|\beta|\sqrt{\beta!}
\|\partial^{\alpha-\beta} f\|_{A,\theta}\|\partial^{\alpha}f\|_{2,\theta+\frac{\gamma}{2}}.
\end{align*}
\end{proof}
Substituting the estimates of $\mathbf{R}_{0}(f),\mathbf{R}_{1}(f)$ into the decomposition \eqref{linear-decomp}, we can end the proof of Proposition \ref{multiply-estimate-weighted}.

\section{Existence and uniqueness of linear Landau equation.}\label{S4}

\begin{proposition}\label{existence}
For $-3< \gamma<0$, $T>0$, $f_0\in L^2(\mathbb{R}^3)$,  the Cauchy problem \eqref{eq-1} admits a unique weak solution
$$
f\in L^{\infty}([0,T]; L^2(\mathbb{R}^3))
$$
satisfying
$$
\|f(t)\|_{L^{2}(\mathbb{R}^3)}\leq 2e^{2C_{0}T}\| f_{0}\|_{L^{2}(\mathbb{R}^3)}
$$
where $C_0$ was defined in Proposition \ref{multiply-estimate-weighted}.
\end{proposition}
\begin{proof}
The existence of the weak solution is similar to that in \cite{AMUXY, MX, Guo-2002}.  We consider the operator
$$
\mathcal{P}=-\partial_t+\mathcal{L}.
$$
For any $\varphi\in C^{\infty}([0,T]; \mathcal{S}(\mathbb{R}^3))$ with $\varphi(T)=0,$ it follows from
\eqref{Remark-1} that
\begin{align*}
(\varphi(t), \mathcal{P}^*\varphi(t))_{L^2(\mathbb{R}^3)}
&=(\partial_{t}\varphi,\varphi)_{L^{2}(\mathbb{R}^{3})}+(\mathcal{L}\varphi,\varphi)_{L^{2}(\mathbb{R}^{3})}\\
&\leq\frac{1}{2}\frac{d}{dt}\|\varphi\|_{L^{2}(\mathbb{R}^{3})}^{2}-\|\varphi\|_{A}^{2}+C_{0}\| \varphi\|_{2,\frac{\gamma}{2}}^{2}.
\end{align*}
Since $\gamma<0$, we have
\begin{align*}
-\frac{1}{2}\frac{d}{dt}\|\varphi\|_{L^{2}(\mathbb{R}^{3})}^{2}+\|\varphi\|_{A}^{2}
&\leq(\varphi(t), \mathcal{P}^*\varphi(t))_{L^2(\mathbb{R}^3)}+C_{0}\| \varphi\|_{L^{2}(\mathbb{R}^{3})}^{2}.
\end{align*}
Which implies that
\begin{align*}
 -\frac{d}{dt}\left(e^{2C_{0}t}\|\varphi\|_{L^{2}(\mathbb{R}^{3})}^{2}\right)+e^{2C_{0}t}\|\varphi\|_{A}^{2}
&\leq 2e^{2C_{0}t}\|\varphi\|_{L^{2}(\mathbb{R}^{3})}\|\mathcal{P}^*\varphi(t)\|_{L^{2}(\mathbb{R}^{3})},
\end{align*}
Then one  can verify that
\begin{equation}\label{injection}
\sup_{0\leq t\le T}\|\varphi(t)\|_{L^{2}(\mathbb{R}^{3})}^{2}
\leq2 e^{2C_0T}\int^T_0\|\mathcal{P}^*\varphi(s)\|_{L^{2}(\mathbb{R}^{3})}ds.
\end{equation}
In the following, we set the vector subspace,
\begin{align*}
&\mathbb{U}=\left\{u=\mathcal{P}\varphi:\varphi\in C^{\infty}([0,T],\mathcal{S}(\mathbb{R}^{6}_{x,v})),\varphi(T)=0\right\}
\subset L^{1}([0,T],L^2(\mathbb{R}^3))).
\end{align*}
Since $f_{0}\in L^2(\mathbb{R}^3)$, we define the linear functional
\begin{align*}
&\mathcal{G}:\mathbb{U}\rightarrow\mathbb{R}
\\&
h=\mathcal{P}^*\varphi\mapsto(f_{0},\varphi(0))_{L^2(\mathbb{R}^3)},
\end{align*}
where $\varphi\in C^{\infty}([0,T],\mathcal{S}(\mathbb{R}^3))$ with $\varphi(T)=0$.
From \eqref{injection}, the operator $\mathcal{P}^*$ is injective. Therefore the linear functional $\mathcal{G}$ is well-defined and moreover, we obtain
\begin{align*}
|\mathcal{\mathcal{G}}(h)|&\leq\|f_0\|_{L^2}\|\varphi(0)\|_{L^2(\mathbb{R}^3)}\leq\|f_0\|_{L^2}
\|\varphi\|_{L^{\infty}([0,T];L^2(\mathbb{R}^3))}
\\&\leq 2e^{2C_0T}\|f_{0}\|_{L^2}\|\mathcal{P}^*\varphi\|_{L^{1}([0,T],L^2(\mathbb{R}^3))}
=2e^{2C_0T}\|f_{0}\|_{L^2}\|h\|_{L^{1}([0,T];L^2(\mathbb{R}^3))}.
\end{align*}
This shows that $\mathcal{G}$ is a continuous linear form on $(\mathbb{U},\|\cdot\|_{L^{1}([0,T],L^2(\mathbb{R}^3))})$.
By using the Hahn-Banach theorem, $\mathcal{G}$ can be extended as a continuous linear form on
$L^{1}([0,T];L^2(\mathbb{R}^3))$ with a norm smaller than $2e^{2C_0T}\|f_{0}\|_{L^2(\mathbb{R}^3)}$.
It follows the Riesz representation theorem  that there exists a unique
 $f\in L^{\infty}([0,T];L^2(\mathbb{R}^3))$
 satisfying
\begin{align*}
\|f\|_{L^{\infty}([0,T];L^2(\mathbb{R}^3))}\leq2e^{2C_0T}\|f_{0}\|_{L^2(\mathbb{R}^3)},
\end{align*}
such that
\begin{align*}
\forall h\in L^{1}([0,T];L^2(\mathbb{R}^3)), \quad \quad \mathcal{G}(h)=\int_{0}^{T}\left(f(t),h(t)\right)_{L^2(\mathbb{R}^3)}dt.
\end{align*}
Which implies that for all $\varphi\in C^{\infty}_{0}((0,T),\mathcal{S}(\mathbb{R}^3))$,
\begin{align*}
&\mathcal{G}(\mathcal{P}^*\varphi)=\int_{0}^{T}\left(f(t),\mathcal{P}^*\varphi(t)\right)_{L^2(\mathbb{R}^3)}dt
=(f_{0},\varphi(0))_{L^2(\mathbb{R}^3)}.
\end{align*}
Therefore, $f\in L^{\infty}([0,T];L^2(\mathbb{R}^3))$ is a weak solution of the Cauchy problem \eqref{eq-1}.
Let $\tilde{f}\in L^{\infty}([0,T];L^2(\mathbb{R}^3))$ be another weak solution of the Cauchy problem \eqref{eq-1} satisfying
\begin{align*}
\|\tilde{f}\|_{L^{\infty}([0,T];L^2(\mathbb{R}^3))}\leq2e^{2C_0T}\|f_{0}\|_{L^2(\mathbb{R}^3)},
\end{align*}
Set $w(t)=f(t)-\tilde{f}(t)$, we have that for all $\varphi\in C^{\infty}_{0}((0,T),\mathcal{S}(\mathbb{R}^3))$,
\begin{align*}
\int_{0}^{T}\left(w(t),\mathcal{P}^*\varphi(t)\right)_{L^2(\mathbb{R}^3)}dt
=0.
\end{align*}
Which shows that $w(t)=0$ in $L^{\infty}([0,T];L^2(\mathbb{R}^3))$.
The proof of Proposition \ref{existence} is completed.

\end{proof}

\section{Analytic smoothing effect for linear Landau equation}\label{S5}

Let $f$ be the solution of Cauchy problem \eqref{eq-1}, that is,
$$
\partial_{t}f=\mathcal{L}f,\quad f|_{t=0}=f_0,
$$
then similar to the estimate of $\mathbf{R}_{0}$ in Proposition \ref{multiply-estimate-weighted} with  $\theta=0$, we have
\begin{align*}
\frac{d}{dt}\|f\|_{L^{2}(\mathbb{R}^{3})}^{2}
&=2(\partial_{t}f,f)_{L^{2}(\mathbb{R}^{3})}=
2(\mathcal{L}f,f)_{L^{2}(\mathbb{R}^{3})}\\
&\leq-2\|f\|_{A}^{2}+2C_{0}\| f\|_{2,\frac{\gamma}{2}}^{2}.
\end{align*}
This implies that (since $\gamma<0$),
$$
\frac{d}{dt}\|f\|_{L^{2}(\mathbb{R}^3)}^{2}+
2\|f\|_{A}^{2}\leq2C_0\| f\|^2_{L^{2}(\mathbb{R}^3)}.
$$
By using the Gronwall inequality, for any $T>0$ and $0<t<T$, we have
\begin{equation}\label{a1+0}
\|f(t)\|_{L^{2}(\mathbb{R}^3)}^{2}+\int_{0}^{t}\| f(s)\|_{A}^{2}ds\leq e^{2C_{0}T}\| f_{0}\|_{L^{2}(\mathbb{R}^3)}^{2}.
\end{equation}
Let $|\alpha|=1$, it follows from Proposition \ref{multiply-estimate-weighted} with  $\theta=\gamma$ that
\begin{equation}\label{a1}
\begin{split}
& \frac{d}{dt}\|t^{\frac{1}{2}}\langle v\rangle^{\frac{\gamma}{2}}\partial^{\alpha} f\|_{L^{2}(\mathbb{R}^3)}^{2}=2\langle\langle v\rangle^{\gamma}\partial^{\alpha}\partial_{t}f,
t\partial^{\alpha} f\rangle+\|\langle v\rangle^{\frac{\gamma}{2}}\partial^{\alpha} f\|_{L^{2}(\mathbb{R}^3)}^{2}\\
&=2\langle\langle v\rangle^{\gamma}\partial^{\alpha} \mathcal{L}f  ,
t\partial^{\alpha} f\rangle_{L^{2}(\mathbb{R}^3)}+\|\langle v\rangle^{\frac{\gamma}{2}}\partial^{\alpha} f\|_{L^{2}(\mathbb{R}^3)}^{2}\\
&\leq-2\|t^{\frac{1}{2}}\partial^{\alpha} f\|_{A,\frac{\gamma}{2}}^{2}+2C_{0}t\| f\|_{A,\frac{\gamma}{2}}\|\partial^{\alpha} f\|_{A,\frac{\gamma}{2}}+\|\langle v\rangle^{\frac{\gamma}{2}}\partial^{\alpha} f\|_{L^{2}(\mathbb{R}^3)}^{2}.
\end{split}
\end{equation}
By using the inequality \eqref{definition2}, one can verify that, for any $0<t<T$ and $\gamma<0$,
\begin{align*}
\| f\|_{A,\frac{\gamma}{2}}\le \| f\|_{A},\quad\|\langle v\rangle^{\frac{\gamma}{2}}\partial^{\alpha} f\|_{L^{2}(\mathbb{R}^3)}^{2}\leq\frac{1}{C_1}\| f\|^2_{A}.
\end{align*}
Substituting back to the estimate \eqref{a1}, we have
\begin{align*}
\frac{d}{dt}\|t^{\frac{1}{2}}\partial^{\alpha} f\|_{2, \frac{\gamma}{2}}^{2}+\|t^{\frac{1}{2}}\partial^{\alpha} f\|_{A,\frac{\gamma}{2}}^{2}\leq \left(C_{0}^{2}t+\frac{1}{C_1}\right)\|f\|_{A}^{2}.
\end{align*}
Integral on $[0, t]$, and using \eqref{a1+0}, one can verify that, for $|\alpha|=1$,
\begin{equation}\label{C+1}
\|{t}^{\frac{1}{2}}\partial^{\alpha} f\|_{2, \frac{\gamma}{2}}^{2}
+\int_{0}^{{t}}\| \tau^{\frac{1}{2}} \partial^{\alpha}u\|_{A}^{2}d\tau
\leq\left(C_{0}^2T+\frac{1}{C_1}\right)e^{2C_{0}T}
\|f_{0}\|_{L^{2}(\mathbb{R}^3)}^{2}\leq C^{4}.
\end{equation}

\begin{proposition}\label{very}
For any $m\in \mathbb{N}$ and $\alpha\in\mathbb{N}^3, |\alpha|=m$, we have, for $0<t\le 1$,
\begin{equation}\label{Assumption2}
\|{t}^{\frac{|\alpha|}{2}}\partial^{\alpha} f\|_{2, \frac{\gamma|\alpha|}{2}}^{2}+\int_{0}^{{t}}\| \tau^{\frac{|\alpha|}{2}}\partial^{\alpha} f\|_{A,\frac{\gamma|\alpha|}{2}}^{2}d\tau\leq {C}^{2|\alpha|+2}(\alpha!)^{2}.
\end{equation}
where $C$ is independent on $\alpha$.
\end{proposition}

This Proposition imply Theorem \ref{trick}.

\begin{proof}
In fact, we have proved that the assumption \eqref{Assumption2} holds true for $m=0,1$ by \eqref{a1+0} and \eqref{C+1}.

Now assume that the assumption \eqref{Assumption2} holds true for $|\alpha|\leq m-1$, that means, for any $|\alpha|\le m-1$, for $0<t\le 1$,
\begin{equation}\label{Assumption2+m+1}
\|{t}^{\frac{|\alpha|}{2}}\partial^{\alpha} f\|_{2, \frac{\gamma|\alpha|}{2}}^{2}+\int_{0}^{{t}}\| \tau^{\frac{|\alpha|}{2}}\partial^{\alpha} f\|_{A,\frac{\gamma|\alpha|}{2}}^{2}d\tau\leq {C}^{2|\alpha|+2}(\alpha!)^{2}.
\end{equation}
We intend to prove the validity of \eqref{Assumption2} for $m$. Firstly
\begin{align*}
&\frac{d}{dt}\|t^{\frac{|\alpha|}{2}}\langle v\rangle^{\frac{\gamma|\alpha|}{2}}\partial^{\alpha} f\|_{L^{2}(\mathbb{R}^3)}^{2}\\
&=2t^{|\alpha|}\langle\langle v\rangle^{\gamma|\alpha|}\partial^{\alpha}\partial_{t}f,
\partial^{\alpha}f\rangle_{L^{2}(\mathbb{R}^3)}+|\alpha|t^{|\alpha|-1}\| \langle v\rangle^{\frac{\gamma|\alpha|}{2}}\partial^{\alpha}  f\|_{L^{2}(\mathbb{R}^3)}^{2}\\
&=2t^{|\alpha|}\langle\langle v\rangle^{\gamma|\alpha|}\partial^{\alpha}\mathcal{L}f,
\partial^{\alpha}f\rangle_{L^{2}(\mathbb{R}^3)}+|\alpha|t^{|\alpha|-1}\| \langle v\rangle^{\frac{\gamma|\alpha|}{2}}\partial^{\alpha}  f\|_{L^{2}(\mathbb{R}^3)}^{2}
\end{align*}
Let  $\theta=\frac{\gamma|\alpha|}{2}$  in Proposition \ref{multiply-estimate-weighted}, it follows that,
$$
\frac{d}{dt}\|t^{\frac{|\alpha|}{2}}\langle v\rangle^{\frac{\gamma|\alpha|}{2}}\partial^{\alpha} f\|_{L^{2}(\mathbb{R}^3)}^{2}
\leq-2\|t^{\frac{|\alpha|}{2}}\partial^{\alpha} f\|_{A,\frac{\gamma|\alpha|}{2}}^{2}+\mathbf{B}(f),
$$
and
\begin{equation}\label{int-1}
\sum_{|\alpha|=m}\|t^{\frac{|\alpha|}{2}}\langle v\rangle^{\frac{\gamma|\alpha|}{2}}\partial^{\alpha} f\|_{L^{2}(\mathbb{R}^3)}^{2}+2\sum_{|\alpha|=m}\int^t_0\|{\tau}^{\frac{|\alpha|}{2}}\partial^{\alpha} f\|_{A,\frac{\gamma|\alpha|}{2}}^{2}d\tau\le \int^t_0 \mathbf{B}(f)d\tau,
\end{equation}
with
\begin{align*}
\mathbf{B}(f)&=|\alpha|t^{|\alpha|-1}\| \langle v\rangle^{\frac{\gamma|\alpha|}{2}}\partial^{\alpha}  f\|_{L^{2}(\mathbb{R}^3)}^{2}\\
&\quad+2C_0
\sum_{\substack{|\beta|\geq1\\
\beta\leq\alpha}}C_{\alpha}^{\beta}\sqrt{\beta!}t^{|\alpha|}
\|\partial^{\alpha-\beta}f\|_{A,\frac{\gamma|\alpha|}{2}}\|\partial^{\alpha}f\|_{A,\frac{\gamma|\alpha|}{2}}\\
&\quad+2C_0\sum_{\substack{|\beta|\geq2\\
\beta\leq\alpha}}C_{\alpha}^{\beta}|\beta|\sqrt{\beta!}t^{|\alpha|}
\|\partial^{\alpha-\beta}f\|_{A,\frac{\gamma|\alpha|}{2}}\|\partial^{\alpha}f\|_{2,\frac{\gamma|\alpha|}{2}+\frac{\gamma}{2}}\\
&\quad+C_0|\gamma||\alpha|\sum_{
\beta\leq\alpha}C_{\alpha}^{\beta}\sqrt{\beta!}t^{|\alpha|}
\|\partial^{\alpha-\beta}f\|_{A,\frac{\gamma|\alpha|}{2}}\|\partial^{\alpha}f\|_{\frac{\gamma|\alpha|}{2}+\frac{\gamma}{2}}\\
&=\mathbf{B}_1(f)+\mathbf{B}_2(f)+\mathbf{B}_3(f)+\mathbf{B}_4(f).
\end{align*}
It follows from  the inequality \eqref{definition2} again,
\begin{align*}
\mathbf{B}_1(f)&=|\alpha|t^{|\alpha|-1}\| \langle v\rangle^{\frac{\gamma|\alpha|}{2}}\partial^{\alpha}  f\|_{L^{2}(\mathbb{R}^3)}^{2}\\
&\le |\alpha| t^{|\alpha-e_{k_0}|}\| \langle v\rangle^{\frac{\gamma|\alpha-e_{k_0}|}{2}}\partial^{\alpha-e_{k_0}}  \nabla f\|_{L^{2}(\mathbb{R}^3)}^{2}\\
&\leq \frac{1}{C_1}|\alpha|\|t^{\frac{|\alpha-e_{k_0}|}{2}}\langle v\rangle^{\frac{\gamma(|\alpha-e_{k_0}|)}{2}}\partial^{\alpha-e_{k_0}}f\|^2_{A},
\end{align*}
where $k_0$ is choose with $\alpha_{k_0}=\max\{\alpha_1, \alpha_2, \alpha_3\}$.
So that, by induction assumption \eqref{Assumption2+m+1} for $|\alpha-e_{k_0}|=m-1$
\begin{align*}
\int^t_0 \mathbf{B}_1(f) d\tau &\le \frac{1}{C_1}|\alpha|\int^t_0\|{\tau}^{\frac{|\alpha-e_{k_0}|}{2}}\langle v\rangle^{\frac{\gamma(|\alpha-e_{k_0}|)}{2}}\partial^{\alpha-e_{k_0}}f\|^2_{A}d\tau\\
&\le \frac{|\alpha|}{C_1}C^{2|\alpha-e_{k_0}|+2}({\alpha-e_k}!)^2
\le \frac{|\alpha|}{C_1(\alpha_k)^2}C^{2|{\alpha}|}({\alpha}!)^2.
\end{align*}
We get then
\begin{equation}\label{Indu-1}
\int^t_0 \mathbf{B}_1(f) d\tau \le \frac{3}{C_1}C^{2|{\alpha}|}({\alpha}!)^2
\end{equation}
For the term $\mathbf{B}_2(f) $, using the fact, for any $\beta\leq\alpha$, $\gamma<0$,
$$
\|\partial^{\alpha-\beta}f\|_{A,\frac{\gamma|\alpha|}{2}}\leq \|\partial^{\alpha-\beta}f\|_{A,\frac{\gamma|\alpha-\beta|}{2}}, $$
by using induction assumption \eqref{Assumption2+m+1} for $|\alpha-\beta|\le m-1$, for $0<t\le 1$,
\begin{align*}
&\int^t_0 \mathbf{B}_2(f) d\tau \le 2C_0 T^m
\sum_{\substack{|\beta|\geq1\\
\beta\leq\alpha}}C_{\alpha}^{\beta}\sqrt{\beta!}\int^t_0
\|\tau^{\frac{|\alpha-\beta|}{2}}\partial^{\alpha-\beta}f\|_{A,\frac{\gamma|\alpha-\beta|}{2}}
\|\tau^{\frac{|\alpha|}{2}}\partial^{\alpha}f\|_{A,\frac{\gamma|\alpha|}{2}}d\tau\\
 &\le 2C_0
\sum_{\substack{|\beta|\geq1\\
\beta\leq\alpha}}C_{\alpha}^{\beta}\sqrt{\beta!}\left(\int^t_0
\|\tau^{\frac{|\alpha-\beta|}{2}}\partial^{\alpha-\beta}f\|^2_{A,\frac{\gamma|\alpha-\beta|}{2}}d\tau\right)^{1/2}
\left(\int^t_0\|\tau^{\frac{|\alpha|}{2}}\partial^{\alpha}f\|^2_{A,\frac{\gamma|\alpha|}{2}}d\tau\right)^{1/2}\\
&\le 4\left(2C_0
\sum_{|\beta|\geq1,
\beta\leq\alpha}C_{\alpha}^{\beta}\sqrt{\beta!}C^{|\alpha-\beta|+1}(\alpha-\beta)! \right)^2+
\frac 12 \int^t_0\|\tau^{\frac{|\alpha|}{2}}\partial^{\alpha}f\|^2_{A,\frac{\gamma|\alpha|}{2}}d\tau.
\end{align*}
We get then, for $0<t\le 1$,
\begin{equation}\label{Indu-2}
\int^t_0 \mathbf{B}_2(f) d\tau \le \tilde{C}_1^2 C^{2|{\alpha}|}({\alpha}!)^2
+
\frac 12 \int^t_0\|\tau^{\frac{|\alpha|}{2}}\partial^{\alpha}f\|^2_{A,\frac{\gamma|\alpha|}{2}}d\tau
\end{equation}
with
\begin{equation}\label{Int-2}
\tilde{C}_1\ge 4C_0
\sum_{|\beta|\geq1}\frac{1}{\sqrt{\beta!}}C^{1-|\beta|}.
\end{equation}
For the term $\mathbf{B}_3(f) $, by using induction assumption \eqref{Assumption2+m+1} for $|\alpha-\beta|\le m-2, |\alpha-e_{k_0}|=m-1$,
\begin{align*}
&\int^t_0 \mathbf{B}_3(f) d\tau \le 2C_0\sum_{\substack{|\beta|\geq2\\
\beta\leq\alpha}}C_{\alpha}^{\beta}|\beta|\sqrt{\beta!}\int^t_0 \tau^{|\alpha|}
\|\partial^{\alpha-\beta}f\|_{A,\frac{\gamma|\alpha|}{2}}\|\partial^{\alpha}f\|_{2,\frac{\gamma|\alpha|}{2}
+\frac{\gamma}{2}}d\tau\\
&\le 2C_0\sum_{\substack{|\beta|\geq2\\
\beta\leq\alpha}}C_{\alpha}^{\beta}|\beta|\sqrt{\beta!}\left(\int^t_0 \|\tau^{\frac{|\alpha-\beta|}{2}}
\partial^{\alpha-\beta}f\|^2_{A,\frac{\gamma|\alpha-\beta|}{2}}d\tau\right)^{1/2}\\
&\qquad\qquad\qquad \qquad\qquad\qquad\times \left(\int^t_0\|\tau^{\frac{|\alpha|-1}{2}}\partial^{\alpha-e_{k_0}}\nabla f\|^2_{2,\frac{\gamma|\alpha|-1}{2}+\frac{\gamma}{2}}d\tau\right)^{1/2}\\
&\le 2C_0\sum_{\substack{|\beta|\geq2\\
\beta\leq\alpha}}C_{\alpha}^{\beta}|\beta|\sqrt{\beta!}\left(\int^t_0 \|\tau^{\frac{|\alpha-\beta|}{2}}
\partial^{\alpha-\beta}f\|^2_{A,\frac{\gamma|\alpha-\beta|}{2}}d\tau\right)^{1/2}\\
&\qquad\qquad\qquad \qquad\qquad\qquad\times \left(\int^t_0\|\tau^{\frac{|\alpha|-1}{2}}\partial^{\alpha-e_{k_0}} f\|^2_{A,\frac{\gamma|\alpha|-1}{2}}d\tau\right)^{1/2}\\
&\le 2C_0\sum_{\substack{|\beta|\geq2\\
\beta\leq\alpha}}C_{\alpha}^{\beta}|\beta|\sqrt{\beta!}C^{|\alpha-\beta|+1}(\alpha-\beta)!\
C^{|\alpha|}(\alpha-e_{k_0})!,
\end{align*}
we get then
\begin{equation}\label{Indu-3}
\int^t_0 \mathbf{B}_3(f) d\tau \le \tilde{C}^2_2 C^{2|{\alpha}|}({\alpha}!)^2
\end{equation}
with
\begin{equation}\label{Int-3}
\tilde{C}_2^2\ge 6C_0
\sum_{|\beta|\geq 2, \beta\le \alpha}\frac{1}{\sqrt{\beta!}}C^{1-|\beta|}\ge 2C_0
\sum_{|\beta|\geq 2, \beta\le \alpha}\frac{|\beta|}{\sqrt{\beta!}\alpha_{k_0}}C^{1-|\beta|}.
\end{equation}
Finally for the term $\mathbf{B}_4(f) $, by using induction assumption \eqref{Assumption2+m+1} for $|\alpha-e_{k_0}|=m-1$,
\begin{align*}
\int^t_0 \mathbf{B}_4(f) d\tau &\le C_0|\gamma||\alpha|\sum_{
0\le \beta\leq\alpha}C_{\alpha}^{\beta}\sqrt{\beta!}\int^t_0\big( \|\tau^{\frac{|\alpha-\beta|}{2}}
\|\partial^{\alpha-\beta}f\|_{A,\frac{\gamma|\alpha|}{2}}\\
&\qquad\qquad\qquad\qquad\qquad \times \|\tau^{\frac{|\alpha-e_{k_0}|}{2}}\partial^{\alpha-e_{k_0}}\nabla
f\|_{\frac{\gamma|\alpha-e_{k_0}|}{2}+\frac{\gamma}{2}}\big)d\tau\\
&\le C_0|\gamma||\alpha|\sum_{
0\le\beta\leq\alpha}C_{\alpha}^{\beta}\sqrt{\beta!}\left(\int^t_0\|\tau^{\frac{|\alpha-\beta|}{2}}
\|\partial^{\alpha-\beta}f\|^2_{A,\frac{\gamma|\alpha-\beta|}{2}}d\tau\right)^{1/2}\\
&\qquad\qquad\qquad\qquad \times \left(\int^t_0\|\tau^{\frac{|\alpha-e_{k_0}|}{2}}\partial^{\alpha-e_{k_0}}
f\|_{A, \frac{\gamma|\alpha-e_{k_0}|}{2}}d\tau\right)^{1/2}\\
&\le C_0|\gamma||\alpha|\left(\int^t_0\|\tau^{\frac{|\alpha|}{2}}
\|\partial^{\alpha}f\|^2_{A,\frac{\gamma|\alpha|}{2}}d\tau\right)^{1/2}C^{|\alpha-e_{k_0}|+1}(\alpha-e_{k_0})!\\
&\qquad+ C_0|\gamma||\alpha|\sum_{
0<\beta\leq\alpha}C_{\alpha}^{\beta}\sqrt{\beta!}C^{|\alpha-\beta|+1}(\alpha-\beta)! C^{|\alpha-e_{k_0}|+1}(\alpha-e_{k_0})!,
\end{align*}
so that
\begin{align*}
\int^t_0 \mathbf{B}_4(f) d\tau &\le \frac{1}{2}\int^t_0\|\tau^{\frac{|\alpha|}{2}}
\|\partial^{\alpha}f\|^2_{A,\frac{\gamma|\alpha|}{2}}d\tau
+4C^2_0|\gamma|^2|\alpha|^2 C^{2|\alpha|}((\alpha-e_{k_0})!)^2\\
&\qquad+ C_0|\gamma||\alpha|\sum_{
0<\beta\leq\alpha}C_{\alpha}^{\beta}\sqrt{\beta!}C^{|\alpha-\beta|+1}(\alpha-\beta)! C^{|\alpha-e_{k_0}|+1}(\alpha-e_{k_0})!,
\end{align*}
We get then, for $0<t\le 1$,
\begin{equation}\label{Indu-4}
\int^t_0 \mathbf{B}_4(f) d\tau \le \tilde{C}_3^2 C^{2|{\alpha}|}({\alpha}!)^2
+
\frac 12 \int^t_0\|\tau^{\frac{|\alpha|}{2}}\partial^{\alpha}f\|^2_{A,\frac{\gamma|\alpha|}{2}}d\tau
\end{equation}
with
\begin{equation}\label{Int-4}
\tilde{C}^2_3\ge 4C^2_0|\gamma|^2|\frac{|\alpha|^2}{\alpha^2_{k_0}}+
\sum_{|\beta|\geq1}\frac{C_0|\gamma||\alpha|}{\sqrt{\beta!}\alpha_{k_0}}C^{1-|\beta|}.
\end{equation}
Taking the constant $C$ satisfying \eqref{C+1}, and
$$
C^2\ge \frac{3}{C_1}+\tilde{C}^2_1+\tilde{C}^2_2+\tilde{C}^2_3,
$$
where the constants are defined by \eqref{Int-2}, \eqref{Int-3} and \eqref{Int-4}.
Combine \eqref{int-1}, \eqref{Indu-1}, \eqref{Indu-2}, \eqref{Indu-3} and \eqref{Indu-4}, we ends the proof of Proposition \ref{very}.
\end{proof}

\section{Appendix}\label{Appendix}
The standard Hermite functions $(\varphi_{n})_{n\in \mathbb{N}}$ are defined for $v\in\mathbb{R}$,
\begin{align*}
\varphi_{n}(v)=\frac{(-1)^n}{\sqrt{2^nn!\sqrt{\pi}}}e^{\frac{v^2}{2}}\frac{d^n}{dv^n}(e^{-\frac{v^2}{2}})
=-\frac{1}{\sqrt{2^nn!\sqrt{\pi}}}(v-\frac{d}{dv})^n(e^{-\frac{v^2}{2}})=\frac{a_+^n\varphi_{0}}{\sqrt{n!}},
\end{align*}
where $a_+$ is the creation operator
$$
a_+=\frac{1}{\sqrt{2}}\Big(v-\frac{d}{dv}\Big).
$$
The family $(\varphi_{n})_{n\in \mathbb{N}}$ is an orthonormal basis of $L^2(\mathbb{R})$. we set for $n\geq0$, $\alpha=(\alpha_1,\alpha_2,\alpha_3)\in \mathbb{N}^3$, $x\in \mathbb{R}$,$v\in \mathbb{R}^3$,\\
  \begin{equation*}
  \psi_{n}(x)=2^{-1/4}\varphi_{n}(2^{-1/2}x),\ \ \ \
 \psi_{n}=\frac{1}{\sqrt{n!}}\left(\frac{x}{2}-\frac{d}{dx}\right)^{n}\psi_{0},
 \end{equation*}
 \begin{equation*}
 \Psi_{\alpha}(v)=\prod_{j=1}^{3}\psi_{\alpha_{j}}(v_j),\ \ \ \
\mathcal{E}_{k}=Span({\Psi_{\alpha}})_{\alpha\in N^3,|\alpha|=k},
 \end{equation*}
with $|\alpha|=\alpha_{1}+\alpha_{2}+\alpha_3$.  The family $(\Psi_{\alpha})_{\alpha\in \mathbb{N}^3}$ is an orthonormal basis of $ L^{2}(\mathbb{R}^3)$ composed by the eigenfunctions of the 3-dimensional harmonic oscillator
 \begin{equation}\label{decomH}
\mathcal{H}=-\Delta_{v}+\frac{|v^{2}|}{4}=\sum_{k\geq0}(k+\frac{3}{2})\mathbb{P}_{k},\quad Id=\sum_{k\geq0}\mathbb{P}_{k},
 \end{equation}
 where $\mathbb{P}_{k}$ stands for the orthogonal projection
$$
\mathbb{P}_kf=\sum_{|\alpha|=k}(f,\Psi_{\alpha})_{L^2(\mathbb{R}_v)}\Psi_{\alpha}.
$$
In particular,
\begin{equation*}
 \Psi_{0}(v)=\frac{1}{(2\pi)^{\frac{3}{4}}}e^{-\frac{|v|^{2}}{4}}=\mu^{1/2}(v),
 \end{equation*}
where $\mu(v)$ is the Maxwellian distribution.    Setting
\begin{equation}\label{H3}
A_{\pm,j}=\frac{v_{j}}{2}\mp\partial_j,\quad 1\leq j\leq 3,
\end{equation}
we have
\begin{equation*}
\Psi_{\alpha}=\frac{1}{\sqrt{\alpha_{1}!\alpha_{2}!\alpha_{3}!}}A^{\alpha_{1}}_{+,1}A^{\alpha_{2}}_{+,2}A^{\alpha_{3}}_{+,3}\Psi_{0},\quad  \alpha=(\alpha_{1},\alpha_{2},\alpha_{3})\in \mathbb{N}^{3},
\end{equation*}
and
\begin{equation}\label{H4}
A_{+,j}\Psi_{\alpha}=\sqrt{\alpha_{j}+1}\Psi_{\alpha+e_{j}},\quad A_{-,j}\Psi_{\alpha}=\sqrt{\alpha_{j}}\Psi_{_{\alpha-e_{j}}}
(=0 \, if\,  \alpha_{j}=0),
\end{equation}
where $(e_{1},e_2,e_{3})$ stands for the canonical basis of $\mathbb{R}^{3}$.  For more details of the Hermite functions, we can refer to \cite{MPX} and the reference theorem.

\section*{Acknowledgments}

The first author is supported by the Fundamental
Research Funds for the Central Universities of China, South-Central University for Nationalities (No. CZT20007). The second author is supported by the NSFC (No.12031006) and the Fundamental
Research Funds for the Central Universities of China.

\end{document}